\documentclass[a4paper,11pt]{article}
\usepackage{graphicx}
\usepackage{multirow}
\usepackage{color}
\usepackage{latexsym}
\usepackage{amssymb}
\usepackage{amsfonts}
\usepackage{amsmath}
\usepackage{amsthm}
\usepackage{indentfirst}
\usepackage{mathrsfs}
\usepackage{tikz}
\usetikzlibrary{patterns}
\usetikzlibrary{hobby,decorations.pathreplacing,arrows.meta}
\usepackage{tikz}
\usepackage{pictex,dcpic}
\usepackage{verbatim}
\usepackage{tikzsymbols}
\usetikzlibrary{patterns}
\usetikzlibrary{positioning}
\usetikzlibrary{fit}

\usepackage[latin1]{inputenc}

\usepackage{fullpage}

\newtheorem {theorem}{Theorem}[section]

\newtheorem {lemma}{Lemma}[section]
\newtheorem {proposition}{Proposition}[section]
\newtheorem {corollary}{Corollary}[section]

\newtheorem {conjecture}{Conjecture}[section]

\theoremstyle{remark}
\newtheorem*{remark}{Remark}

\date{}

\author{Matteo Cervetti
\thanks{Dipartimento di Matematica, Universit\`a di Trento, Via Sommarive 14, 38123 Povo (Trento), Italy.
{\tt \ matteo.cervetti@unitn.it}}
\and
Luca Ferrari
\thanks{Dipartimento di Matematica e Informatica ``U. Dini'',
Universit\`a degli Studi di Firenze, Viale
G.B. Morgagni 65, 50134 Firenze, Italy. {
\tt\ luca.ferrari@unifi.it} Member of the INdAM research group GNCS; partially supported by the 2019 INdAM-GNCS project "Studio di propriet\'a combinatoriche di linguaggi formali ispirate dalla biologia e da strutture bidimensionali".}}

\title{Pattern avoidance in the matching pattern poset}

\begin{document}

\maketitle

\begin{abstract}
A matching of the set $[2n]=\{1,2,...,2n\}$ is a partition of $[2n]$ into blocks with two elements, i.e. a graph on $[2n]$ such that every vertex has degree one. Given two matchings $\sigma$ and $\tau$, we say that $\sigma$ is a $\emph{pattern}$ of $\tau$ when $\sigma$ can be obtained from $\tau$ by deleting some of its edges and consistently relabelling the remaining  vertices. This is a partial order relation turning the set of all matchings into a poset, which will be called the \emph{matching pattern poset}. In this paper, we continue the study of classes of pattern avoiding matchings, initiated by  Chen, Deng, Du, Stanley and Yan (2007), Jelinek and Mansour (2010), Bloom and Elizalde (2012). In particular, we work out explicit formulas to enumerate the class of matchings avoiding two new patterns, obtained by juxtaposition of smaller patterns, and we describe a recursive formula for the generating function of the class of matchings avoiding the lifting of a pattern and two additional patterns. Finally, we introduce the notion of unlabeled pattern, as a combinatorial way to collect patterns, and we provide enumerative formulas for two classes of matchings avoiding an unlabeled pattern of order three. In one case, the enumeration follows from an interesting bijection between the matchings of the class and ternary trees. 
\end{abstract}

\section{Introduction}\label{Intro}

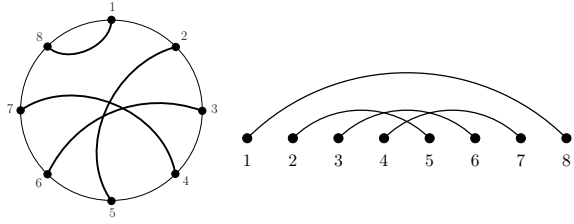
\begin{figure}\begin{center}\scalebox{0.6}{
\begin{tikzpicture}

\node at (0,0) {\scalebox{0.4}{
\begin{tikzpicture}[x=1.00mm, y=1.00mm, inner xsep=0pt, inner ysep=0pt, outer xsep=0pt, outer ysep=0pt]
\path[line width=0mm] (70.51,-1.78) rectangle +(118.56,122.55);
\definecolor{L}{rgb}{0,0,0}
\path[line width=0.30mm, draw=L] (129.69,60.06) circle (49.81mm);
\definecolor{F}{rgb}{0,0,0}
\path[line width=0.30mm, draw=L, fill=F] (79.89,59.98) circle (2.00mm); 
\path[line width=0.30mm, draw=L, fill=F] (129.69,110.01) circle (2.00mm); 
\path[line width=0.30mm, draw=L, fill=F] (179.50,59.74) circle (2.00mm); 
\path[line width=0.30mm, draw=L, fill=F] (129.69,9.94) circle (2.00mm); 
\path[line width=0.30mm, draw=L, fill=F] (164.74,24.93) circle (2.00mm); 
\path[line width=0.30mm, draw=L, fill=F] (164.97,95.02) circle (2.00mm); 
\path[line width=0.30mm, draw=L, fill=F] (94.65,24.70) circle (2.00mm); 
\path[line width=0.30mm, draw=L, fill=F] (94.65,95.25) circle (2.00mm); 
\draw(128.54,114.85) node[anchor=base west]{\fontsize{19.23}{17.07}\selectfont 1};
\draw(167.97,98.71) node[anchor=base west]{\fontsize{19.23}{17.07}\selectfont 2};
\draw(184.57,58.59) node[anchor=base west]{\fontsize{19.23}{17.07}\selectfont 3};
\draw(168.89,18.93) node[anchor=base west]{\fontsize{19.23}{17.07}\selectfont 4};
\draw(128.54,1.30) node[anchor=base west]{\fontsize{19.23}{17.07}\selectfont 5};
\draw(88.19,17.67) node[anchor=base west]{\fontsize{19.23}{17.07}\selectfont 6};
\draw(72.51,58.71) node[anchor=base west]{\fontsize{19.23}{17.07}\selectfont 7};
\draw(88.65,98.37) node[anchor=base west]{\fontsize{19.23}{17.07}\selectfont 8};
\draw[line width=3pt] (129.69,110.01) to [out=270,in=-45] (94.65,95.25) ;
\draw[line width=3pt] (164.97,95.02) to [out=225-30,in=45+90-10] (129.69,9.94) ;
\draw[line width=3pt] (179.50,59.74) to [out=270-90-20,in=-45+90+20] (94.65,24.70);
\draw[line width=3pt] (164.74,24.93) to [out=270-90-90+10,in=-45+90-10] (79.89,59.98) ;
\end{tikzpicture}} };

\node at (6.5,0) {
\begin{tikzpicture}
\draw [fill] (1,0) circle [radius=0.1];
\draw [fill] (2,0) circle [radius=0.1];
\draw [fill] (3,0) circle [radius=0.1];
\draw [fill] (4,0) circle [radius=0.1];
\draw [fill] (5,0) circle [radius=0.1];
\draw [fill] (6,0) circle [radius=0.1];
\draw [fill] (7,0) circle [radius=0.1];
\draw [fill] (8,0) circle [radius=0.1];
\node at (1,-0.5) {1};
\node at (2,-0.5) {2};
\node at (3,-0.5) {3};
\node at (4,-0.5) {4};
\node at (5,-0.5) {5};
\node at (6,-0.5) {6};
\node at (7,-0.5) {7};
\node at (8,-0.5) {8};
\draw[thick] (2,0) to [out=45,in=135] (5,0);
\draw[thick] (3,0) to [out=45,in=135] (6,0);
\draw[thick] (4,0) to [out=45,in=135] (7,0);
\draw[thick] (1,0) to [out=45,in=135] (8,0);
\end{tikzpicture}};

\end{tikzpicture}}\end{center}
\caption{The circular chord diagram representing the perfect matching $\{\{1,8\},\{2,5\},\{3,6\},\{4,7\}\}$ on the set $[8]$ and the corresponding linear chord diagram.}\label{A}
\end{figure}

\begin{figure}\begin{center}\scalebox{0.8}{
\begin{tikzpicture}
\draw [fill] (1,0) circle [radius=0.1];
\draw [fill] (2,0) circle [radius=0.1];
\draw [fill] (3,0) circle [radius=0.1];
\draw [fill] (4,0) circle [radius=0.1];
\draw [fill] (5,0) circle [radius=0.1];
\draw [fill] (6,0) circle [radius=0.1];
\draw [fill] (7,0) circle [radius=0.1];
\draw [fill] (8,0) circle [radius=0.1];
\node at (1,-0.5) {1};
\node at (2,-0.5) {2};
\node at (3,-0.5) {1};
\node at (4,-0.5) {2};
\node at (5,-0.5) {3};
\node at (6,-0.5) {4};
\node at (7,-0.5) {3};
\node at (8,-0.5) {4};
\draw[thick] (1,0) to [out=45,in=135] (3,0);
\draw[thick] (2,0) to [out=45,in=135] (4,0);
\draw[thick] (5,0) to [out=45,in=135] (7,0);
\draw[thick] (6,0) to [out=45,in=135] (8,0);
\end{tikzpicture}}\end{center}
\caption{Encoding the matching $\{\{1,3\},\{2,4\},\{5,7\},\{6,8\}\}$ in the sequence $12123434$.}\label{sequence}
\end{figure}
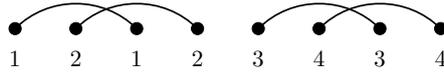

Let $n\in \mathbb{N}^{*}$ and set as usual $[n]=\{1,2,...,n\}$. A $\emph{matching}$ of $[2n]$ is a partition of $[2n]$ into blocks having two elements.  Note that a matching of $[2n]$ is the same as a graph on $[2n]$ such that every vertex has degree one, hence we will borrow some standard terminology from graph theory, as well as the usual representation of graphs using diagrams consisting of dots and lines. In particular, every matching of $[2n]$ will be represented either by a circular or by a linear chord diagram, as shown in Figure $\ref{A}$.
 Let $\tau$ be a matching of $[2n]$. The integer $n$, i.e. the number of edges of $\tau$, will be called the $\emph{order}$ of $\tau$ and will be denoted by $|\tau|$. The set of all matchings will be denoted by $\cal{M}$ and the set of all matchings of order $n$ will be denoted by $\mathcal{M}_{n}$. Given $e\in \tau$, the integers $\min(e)$ and $\max(e)$ will be called the $\emph{left}$ $\emph{vertex}$ and the $\emph{right}$ $\emph{vertex}$ of $e$ respectively. Given a subset $S$ of $\tau$ and $e\in S$, we will say that $e$ is the $\emph{leftmost}$ (respectively $\emph{rightmost}$) edge of $S$ when $\min(e)\leq \min(f)$ (respectively $\max(e)\geq \max(f)$) for every $f\in S$. Following $\cite{JM}$, we will represent $\tau$  by means of the unique integer sequence $\tilde{\tau}\in [n]^{2n}$ such that $\tilde{\tau}_{\min(e)}=\tilde{\tau}_{\max(e)}$ and $\tilde{\tau}_{\min(e)}<\tilde{\tau}_{\min(f)}$ for every $e,f\in \tau$ such that $\min(e)<\min(f)$. Using this encoding, the vertices of $\tau$ are represented by the elements of $\tilde{\tau}$ and two vertices of $\tau$ are connected by an edge when the corresponding components of $\tilde{\tau}$ are equal (see Figure $\ref{sequence}$).
In the following, we will always identify matchings with their corresponding integer sequences. Let $\sigma$ and $\tau$ be matchings. The matching $\sigma(\tau+|\sigma|)$ will be called the $\emph{juxtaposition}$ of $\sigma$ and $\tau$ (where $\tau+|\sigma|$ denotes the sequence obtained from $\tau$ by adding $|\sigma|$ to each of its elements). Its linear chord diagram  can be indeed represented by juxtaposing the linear chord diagrams of $\sigma$ and $\tau$, respectively. The matching $1(\tau+1)1$ will be called the $\emph{lifting}$ of $\tau$. Its linear chord diagram can be represented by nesting the linear chord diagram of $\sigma$ into an additional edge. The matching obtained from the sequence $\tau_{n}...\tau_{2}\tau_{1}$ by suitably renaming its elements so to obtain a valid matching will be called the $\emph{reversal}$ of $\tau$ and denoted by $\overline{\tau}$.  Its linear chord diagram can be represented by reflecting the linear chord diagram of $\tau$ along a vertical line.

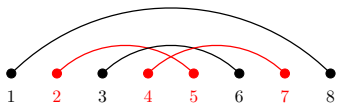
\begin{figure}\begin{center}\scalebox{0.6}{
\begin{tikzpicture}
\draw [fill] (1,0) circle [radius=0.1];
\draw [fill,red] (2,0) circle [radius=0.1];
\draw [fill] (3,0) circle [radius=0.1];
\draw [fill,red] (4,0) circle [radius=0.1];
\draw [fill,red] (5,0) circle [radius=0.1];
\draw [fill] (6,0) circle [radius=0.1];
\draw [fill,red] (7,0) circle [radius=0.1];
\draw [fill] (8,0) circle [radius=0.1];
\node at (1,-0.5) {1};
\node [red] at (2,-0.5) {2};
\node at (3,-0.5) {3};
\node [red] at (4,-0.5) {4};
\node [red] at (5,-0.5) {5};
\node at (6,-0.5) {6};
\node [red] at (7,-0.5) {7};
\node at (8,-0.5) {8};
\draw[thick,red] (2,0) to [out=45,in=135] (5,0);
\draw[thick] (3,0) to [out=45,in=135] (6,0);
\draw[thick,red] (4,0) to [out=45,in=135] (7,0);
\draw[thick] (1,0) to [out=45,in=135] (8,0);
\end{tikzpicture}}\end{center}
\caption{(2,4,5,7) is an occurrence of the perfect matching $\{\{1,3\},\{2,4\}\}$ in the perfect matching $\{\{1,8\},\{2,5\},\{3,6\},\{4,7\}\}$.}\label{MP}
\end{figure}

    Given $k\in \mathbb{N}^{*}$, let $\sigma$ be a matching of $[2k]$ and $i=(i_{1},...,i_{k})\in [2n]^{2k}$. We say that $i$ is an $\emph{occurrence}$ of $\sigma$ in $\tau$ when $i_{1}<i_{2}<...<i_{2k}$ and $\{i_{p},i_{q}\}\in \tau$ if and only if $\{p,q\}\in \sigma$, for every $p,q\in [2k]$ (see Figure $\ref{MP}$).
  We say that $\sigma$ is a $\emph{pattern}$ of $\tau$, and write $\sigma\leq \tau$, when there is an occurrence of $\sigma$ in $\tau$, and that $\tau$ $\emph{avoids}$ $\sigma$ otherwise. The relation $\leq$ is a partial order turning the set of all matchings into a poset, which we will call the $\emph{matching}$ $\emph{pattern}$ $\emph{poset}$. If $S$ is a set of matchings, the class of all matchings avoiding every pattern in $S$ will be denoted by $\mathcal{M}(S)$, the class of all matchings in $\mathcal{M}(S)$ of order $n$ will be denoted by $\mathcal{M}_{n}(S)$ and the generating function of the sequence $\{|\mathcal{M}_{n}(S)|\}_{n\in \mathbb{N}}$ will be denoted by $\mathcal{M}(S,z)$.  We say that $\sigma$ and $\tau$ are $\emph{Wilf-equivalent}$ when $\mathcal{M}(\sigma,z)=\mathcal{M}(\tau,z)$.

In Section $\ref{PA}$  we investigate  classes of the form $\mathcal{M}(\sigma(\tau+|\sigma|))$ and $\mathcal{M}(1(\sigma+1)1,\chi,\overline{\chi})$, providing a general approach which yields enumerative formulas for some patterns $\sigma,\tau$ and $\chi$. Following a recursive approach already described in $\cite{JM}$, we reduce the enumeration of $\mathcal{M}(\sigma(\tau+|\sigma|))$ to the enumeration of a specific class of matchings $\mu(\sigma)$ (depending on $\sigma$) and the class $\mathcal{M}(\tau)$, finding an explicit answer for the prefix $\sigma=1212$. Moreover, we introduce a suitable pattern $\chi=123132$ to relate the generating function of $\mathcal{M}(1(\sigma+1)1,\chi,\overline{\chi})$ to the generating function of $\mathcal{M}(\sigma,\chi,\overline{\chi})$.

In Section $\ref{UPA}$ we introduce the notion of \emph{unlabeled matching}, which is an equivalence class of matchings having the same unlabeled circular chord diagram. This seems a reasonable and combinatorially meaningful way to collect patterns. As a first result concerning unlabeled pattern avoidance, we provide enumerative formulas for two classes of matchings avoiding an unlabeled pattern of order three, as well as a bijection between matchings avoiding a certain unlabeled pattern and ternary trees.

Finally, Section $\ref{Int}$ is a first step towards the study of the enumerative combinatorics of the intervals in the matching pattern poset, and Section $\ref{CFW}$ provides some hints for further work.

\section{Previous work} \label{Work}

Given a permutation $\sigma$ of $[n]$, we can construct a matching of $[2n]$ by connecting the vertices $\{1,...,n\}$ with the vertices $\{n+1,...,2n\}$ in the order prescribed by $\sigma$, thus obtaining the matching corresponding to the integer sequence $12...n\sigma_{1}...\sigma_{n}$. A matching of this kind will be called a $\emph{permutational}$ $\emph{matching}$ and it is immediate to notice that a matching of $[2n]$ is permutational if and only if it avoids the pattern $1122$, so that $|\mathcal{M}_{n}(1122)|=n!$. Two remarkable examples of permutational matchings are $123...n123...n$ and $123...nn...321$, which will be called the $\emph{totally}$ $\emph{crossing}$  and the $\emph{totally}$ $\emph{nesting}$ matching of $[2n]$ respectively. It is easy to see that the map sending every permutation to the corresponding permutational matching is a poset embedding, hence we can regard the permutation pattern poset as a subposet of the matching pattern poset. Throughout this paper we will denote by $C_{n}=\frac{1}{n+1}\binom{2n}{n}$ the $n^{th}$ Catalan number  and  by $C(z)$ the generating function of Catalan numbers. As for other enumerative results on pattern avoidance, it is well known that noncrossing matchings have a Catalan structure, therefore $|\mathcal{M}_{n}(1212)|=C_{n}$, and it is also well known that nonnesting matchings are counted by the same sequence, that is $|\mathcal{M}_{n}(1221)|=C_{n}$.   More surprisingly, it was proved in $\cite{CDDSY}$ that the matchings $123...k123...k$ and $123...kk...321$ are Wilf-equivalent for every $k\in \mathbb{N}^{*}$.  No closed formula for the number of matchings avoiding these patterns is available in general, although it was proved in $\cite{GB}$ that $|\mathcal{M}_{n}(123123)|=C_{n}C_{n+2}-C_{n+1}^{2}$.  Furthermore, Wilf-equivalences between several classes of patterns are established in $\cite{JM}$ through bijective methods; for instance, as an immediate consequence of Lemmas 3.7 and 3.10 in that paper, one can  deduce the following useful fact.

\begin{proposition}\label{W} If $\sigma$ and $\sigma'$ are two Wilf-equivalent matchings and $\tau$ and $\tau'$ are two Wilf-equivalent matchings, then $\sigma(\tau+|\sigma|)$ and $\sigma'(\tau'+|\sigma'|)$ are Wilf-equivalent.
\end{proposition}

Moreover, the same paper also contains an enumerative result which reduces the enumeration of $\mathcal{M}(11(\sigma+1))$ to the enumeration of $\mathcal{M}(\sigma)$ in a recursive fashion. Finally, aside for the classes of patterns mentioned above, the only (up to Wilf-equivalence) further class of matchings avoiding a small pattern has been enumerated in $\cite{BE}$, proving that
$$\mathcal{M}(123132,z)=\frac{54z}{1+36z-(1-12z)^{\frac{3}{2}}}.$$
In the same paper, some enumerative results are given for most of the classes of matchings avoiding a pair of permutational patterns of order three.
Nevertheless, enumerating all the remaining classes of matchings avoiding a single patterns of  order three remains an open problem and it is likely to be a hard one. Indeed, it is suggested in $\cite{BE}$ that enumeration of the class of matchings avoiding the pattern $123231$ could be related to the enumeration of the class of permutations avoiding $1324$, which is considered to be a very hard problem.

\section{Pattern avoidance}\label{PA}

\subsection{The juxtaposition of two patterns}

Let $\sigma$ and $\tau$ be two matchings. In this section we investigate the class of matchings avoiding the juxtaposition of $\sigma$ and $\tau$. To this purpose, we define a set of matchings depending on $\sigma$.  Let $n\in \mathbb{N}$ and $\lambda$ be a matching of order $n$. We will say that $\lambda$ \emph{minimally} \emph{contains} $\sigma$ when it contains $\sigma$ and the matching obtained from $\lambda$ by deleting its rightmost edge does not contain $\sigma$. Denote by $\mu(\sigma)$ the set of matchings minimally containing $\sigma$, by $\mu_{n}(\sigma)$ the set of elements in $\mu(\sigma)$ with order $n$ and by $\mu(\sigma,z)$ the corresponding generating function. Generalizing an approach already used in $\cite{JM}$, the following formula allows us to relate the enumeration of $\mathcal{M}(\sigma(\tau+|\sigma|))$ to the enumeration of $\mathcal{M}(\sigma)$, $\mu(\sigma)$ and $\mathcal{M}(\tau)$.

\begin{proposition}\label{J} Let $\sigma$ and $\tau$ be matchings and $n\in \mathbb{N}$, with $n\geq |\sigma|$. Then
\begin{equation}\label{EJ}
|\mathcal{M}_{n}(\sigma(\tau+|\sigma|))|=|\mathcal{M}_{n}(\sigma)|+\sum_{\ell=|\sigma|}^{n}\sum_{k=0}^{n-\ell}\binom{2\ell+k-1}{k}\binom{2n-2\ell-k}{k}k!|\mu_{\ell}(\sigma)||\mathcal{M}_{n-\ell-k}(\tau)|
\end{equation}
\end{proposition}

\proof Given $\lambda \in \mathcal{M}_{n}(\sigma(\tau+|\sigma|))$, then either $\lambda\in \mathcal{M}_{n}(\sigma)$ or $\sigma$ is a pattern of $\lambda$. From now on we assume that the latter case occurs, since the former one is taken into account by the first summand in the right hand side of $(\ref{EJ})$. For $h\in [2n]$, we denote by $\lambda_{\leq h}$ the pattern of $\lambda$ consisting of all the edges of $\lambda$ with both vertices smaller than or equal to $h$, by $\lambda_{\geq h}$ the pattern of $\lambda$ consisting of all the edges of $\lambda$ with both vertices bigger than or equal to $h$ and finally by $\lambda_{h}$ the pattern of $\lambda$ consisting of all the edges of $\lambda$ that are neither in $\lambda_{\leq h}$ nor in $\lambda_{\geq h}$. Note that an edge of $\lambda$ belongs to $\lambda_{h}$ if and only if its left vertex is smaller than or equal to $h$ and its right vertex is bigger than or equal to $h$.

\begin{center}
\begin{tikzpicture}
\draw [fill] (1,0) circle [radius=0.1];
\draw [fill] (2,0) circle [radius=0.1];
\draw [fill] (4,0) circle [radius=0.1];
\draw [fill] (4,0) circle [radius=0.1];
\draw [fill] (5,0) circle [radius=0.1];
\draw [fill] (6,0) circle [radius=0.1];
\draw [fill] (7,0) circle [radius=0.1];
\draw [fill] (10,0) circle [radius=0.1];
\draw [fill] (11,0) circle [radius=0.1];
\draw[dashed] (1,0) to [out=45,in=135] (11,0) -- (6,0) to [out=135,in=45] (4,0) -- (1,0) ;
\draw[dashed, fill=gray!25, pattern=crosshatch dots] (7,0) to [out=50,in=130] (10,0) -- (7,0);
\draw[dashed, fill=gray!25, pattern=crosshatch dots] (2,0) to [out=50,in=130] (5,0) -- (2,0);
\node at (3.5,0.3) {$\lambda_{\leq h}$};
\node at (8.5,0.3) {$\lambda_{\geq h}$};
\node at (6,1) {$\lambda_{h}$};
\end{tikzpicture}
\end{center}
 Now let $h$ denote the smallest integer such that $\lambda_{\leq h}$ contains an occurrence of $\sigma$ (of course such an $h$ exists because $\lambda_{\leq 2n}=\lambda$) and let $\ell$ be the  order of $\lambda_{\leq h}$. Then, by definition, $\lambda_{\leq h}\in \mu_{\ell}(\sigma)$ and $\ell\in \{|\sigma|,...,n\}$, hence there are $|\mu_{\ell}(\sigma)|$ possible choices for $\lambda_{\leq h}$. Furthermore, $\lambda_{h}\in \mathcal{M}_{k}(1122)$ for some $k\in \{0,...,n-\ell\}$, hence there are $|\mathcal{M}_{k}(1122)|=k!$ possible choices for $\lambda_{h}$. Moreover, $\lambda_{\geq h}\in \mathcal{M}_{n-\ell-k}(\tau)$ because $\lambda\in \mathcal{M}_{n}(\sigma(\tau+|\sigma|))$, hence there are $|\mathcal{M}_{n-\ell-k}(\tau)|$ possibile choices for $\lambda_{\geq h}$. Finally, notice that $h=2\ell+k$ and that the vertex $h$ necessarily belongs to $\lambda_{\leq h}$, hence the left vertices of the edges of $\lambda_{h}$ can be chosen among the vertices of $\lambda$ smaller than $h$ in $\binom{2\ell+k-1}{k}$ ways. Similarly,  the right vertices of the edges of $\lambda_{h}$ can be chosen among the vertices of $\lambda$ bigger than $h$ in $\binom{2n-2\ell-k}{k}$ ways. This explains the factors in the remaining summands of the right hand side of $(\ref{EJ})$ and concludes the proof.
\endproof

Unfortunately, Formula $(\ref{EJ})$ is not very informative, as enumerating $\mu(\sigma)$ is often  as difficult as enumerating $\mathcal{M}(\sigma)$ itself. Nevertheless, one might still hope that this task can be achieved for some special prefixes $\sigma$. For instance,   note that, for $\sigma=11$ and $n\in \mathbb{N}^{*}$, we easily recover the formula
$$|\mathcal{M}_{n}(11(\tau+1))|=\sum_{k=1}^{n}k!\binom{2n-k-1}{k-1}|\mathcal{M}_{n-k}(\tau)|$$
which can be found in $\cite{JM}$. The next proposition shows that the prefix $\sigma=1212$ can be also succesfully addressed.

\begin{proposition}{\label{1212}} Let $n\in \mathbb{N}$, with $n\geq 2$, then
\begin{itemize}
\item[(i)] \begin{equation}\label{RF}|\mu_{n}(1212)|=\sum_{k=0}^{n-2}(2k+1)C_{k}C_{n-k-2}+C_{k}|\mu_{n-k-1}(1212)|;\end{equation}
\item[(ii)] \begin{equation}\label{GF}\mu(1212,z)=\frac{C(z)-1}{(1-2zC(z))(1-zC(z))};\end{equation}
\item[(iii)] $|\mu_{n}(1212)|$ is the $(n-1)^{th}$ term of sequence $\mathrm{A002054}$ in $\cite{S}$, i.e.
$$|\mu_{n}(1212)|=\binom{2n-1}{n-2}.$$
\end{itemize}
\end{proposition}

\proof

 $(i)$  Let $\lambda\in \mu_{n}(1212)$ and let $\hat{\lambda}$ denote the matching obtained from $\lambda$ by removing its rightmost edge, so that $\hat{\lambda}\in \mathcal{M}_{n-1}(1212)$. Using the standard decomposition of noncrossing matchings, we can write $\hat{\lambda}=\pi (|\pi|+1)(\sigma+|\pi|+1)(|\pi|+1)$, where $\pi\in \mathcal{M}_{n-k-2}(1212)$ and $\sigma\in \mathcal{M}_{k}(1212)$ for some $k\in \{0,...,n-2\}$, as shown in the picture below:
\begin{center}
\begin{tikzpicture}
\draw [fill] (1,0) circle [radius=0.1];
\draw [fill] (5,0) circle [radius=0.1];
\draw [fill] (6,0) circle [radius=0.1];
\draw [fill] (6.5,0) circle [radius=0.1];
\draw [fill] (9.5,0) circle [radius=0.1];
\draw [fill] (10,0) circle [radius=0.1];
\draw[thick] (6,0) to [out=45,in=135] (10,0);
\draw[dashed, fill=gray!25, pattern=crosshatch dots] (6.5,0) to [out=45,in=135] (9.5,0) -- (6.5,0);
\draw[dashed, fill=gray!25, pattern=crosshatch dots] (1,0) to [out=45,in=135] (5,0) -- (1,0);
\node at (3,0.3) {$\pi$};
\node at (8,0.3) {$\sigma$};
\end{tikzpicture}
\end{center}
Since $1212\leq \lambda$, there are two cases:
\begin{itemize}
\item[$\bullet$] The rightmost edge of $\lambda$ crosses the rightmost edge of $\hat{\lambda}$, as shown in the picture below:
 \begin{center}
\begin{tikzpicture}
\draw [fill] (1,0) circle [radius=0.1];
\draw [fill] (5,0) circle [radius=0.1];
\draw [fill] (6,0) circle [radius=0.1];
\draw [fill] (6.5,0) circle [radius=0.1];
\draw [fill] (8,0) circle [radius=0.1];
\draw [fill] (9.5,0) circle [radius=0.1];
\draw [fill] (10,0) circle [radius=0.1];
\draw [fill] (11,0) circle [radius=0.1];
\draw[thick] (6,0) to [out=45,in=135] (10,0);
\draw[thick] (8,0) to [out=50,in=130] (11,0);
\draw[dashed, fill=gray!25, pattern=crosshatch dots] (6.5,0) to [out=45,in=135] (9.5,0) -- (6.5,0);
\draw[dashed, fill=gray!25, pattern=crosshatch dots] (1,0) to [out=45,in=135] (5,0) -- (1,0);
\node at (3,0.3) {$\pi$};
\node at (8,0.3) {$\sigma$};
\end{tikzpicture}
\end{center}
In this case the left vertex of the rightmost edge of $\lambda$ can be inserted in all the possible places between the vertices of the rightmost edge of $\hat{\lambda}$, that are $2k+1$ possible places, and therefore there are $(2k+1)C_{k}C_{n-k-2}$ possible choices for $\lambda$.

\item[$\bullet$] The rightmost edge of $\lambda$ crosses an edge of $\pi$, as shown in the picture below:
 \begin{center}
\begin{tikzpicture}
\draw [fill] (1,0) circle [radius=0.1];
\draw [fill] (5,0) circle [radius=0.1];
\draw [fill] (6,0) circle [radius=0.1];
\draw [fill] (6.5,0) circle [radius=0.1];
\draw [fill] (3,0) circle [radius=0.1];
\draw [fill] (9.5,0) circle [radius=0.1];
\draw [fill] (10,0) circle [radius=0.1];
\draw [fill] (11,0) circle [radius=0.1];
\draw[thick] (6,0) to [out=45,in=135] (10,0);
\draw[thick] (3,0) to [out=50,in=130] (11,0);
\draw[dashed, fill=gray!25, pattern=crosshatch dots] (6.5,0) to [out=45,in=135] (9.5,0) -- (6.5,0);
\draw[dashed, fill=gray!25, pattern=crosshatch dots] (1,0) to [out=45,in=135] (5,0) -- (1,0);
\node at (3,0.3) {$\pi$};
\node at (8,0.3) {$\sigma$};
\end{tikzpicture}
\end{center}
In this case, the pattern of $\lambda$ consisting of $\pi$ and the rightmost edge of $\lambda$ minimally contains $1212$, hence there are $C_{k}|\mu_{n-k-1}(1212)|$ possible choices  for $\lambda$.

\end{itemize}
Now summing over $k\in \{0,...,n-2\}$ we find $(\ref{RF})$.

  $(ii)$  It follows from $(i)$ that $$\mu(1212,z)=\sum_{n\geq 2}\left(\sum_{k=0}^{n-2}(2k+1)C_{k}C_{n-k-2}\right)z^{n}+\sum_{n\geq 2}\left(\sum_{k=0}^{n-2}C_{k}|\mu_{n-k-1}(1212)|\right)z^{n}$$
Now
\begin{align*}
\sum_{n\geq 2}\left(\sum_{k=0}^{n-2}(2k+1)C_{k}C_{n-k-2}\right)z^{n}
 & =z^{2}\sum_{n\geq 0}\left(\sum_{k=0}^{n}(2k+1)C_{k}C_{n-k}\right)z^{n}\\
 & =2z^{2}\sum_{n\geq 0}\left(\sum_{k=0}^{n}kC_{k}C_{n-k}\right)z^{n}+ z^{2}\sum_{n\geq 0}\left(\sum_{k=0}^{n}C_{k}C_{n-k}\right)z^{n}\\
 & =z^{2}(2zC'(z)C(z)+C(z)^{2})\\
\end{align*}
On the other hand, we have $C(z)=1+zC(z)^{2}$, therefore
$$2zC(z)C'(z)+C(z)^{2}=(zC(z)^{2})'=(C(z)-1)'=C'(z)$$
hence
$$C'(z)=\frac{C(z)^{2}}{1-2zC(z)}$$
and finally
$$\sum_{n\geq 2}\left(\sum_{k=0}^{n-2}(2k+1)C_{k}C_{n-k-2}\right)z^{n}=\frac{z^{2}C(z)^{2}}{1-2C(z)}=\frac{z(C(z)-1)}{1-2zC(z)}.$$ 
Similarly, we get
\begin{align*}
\sum_{n\geq 2}\left(\sum_{k=0}^{n-2}C_{k}|\mu_{n-k-1}(1212)|\right)z^{n}
 & =z^{2}\sum_{n\geq 0}\left(\sum_{k=0}^{n}C_{k}|\mu_{n-k+1}(1212)|\right)z^{n} \\
 & =z^{2}C(z) \sum_{n\geq 0}|\mu_{n+1}(1212)|z^{n}\\
 & =zC(z)\mu(1212,z)\\
\end{align*}
Summing up, we have
$$\mu(1212,z)=\frac{z(C(z)-1)}{1-2zC(z)}+zC(z)\mu(1212,z),$$
hence
$$\mu(1212,z)=\frac{z(C(z)-1)}{(1-2zC(z))(1-zC(z))}.$$
$(iii)$ The generating function for sequence A002054 can be found in $\cite{S}$ and is given by $$f(z)=\frac{zC(z)^{3}}{1-2zC(z)}.$$
On the other hand
$$zC(z)^{3}(1-zC(z))=C(z)(C(z)-1)(1-zC(z))=$$ $$C(z)(C(z)-1-zC(z)^{2}+zC(z))=zC(z)^{2}=C(z)-1,$$ hence $\mu(1212,z)=zf(z)$, thus proving $(\ref{GF})$.
\endproof

Unfortunately, we have not been able to provide a neat combinatorial argument to explain the appearance of the binomial coefficient in Proposition $\ref{1212}$. However, observe that, as a byproduct, we also find the following identity:
$$\sum_{k=1}^{n-1}\sum_{i=1}^{k}\sum_{\substack{\alpha\in (\mathbb{N}^{*})^{k}\\ |\alpha|=n}}(2\alpha_{i}-1)C_{\alpha_{1}-1}...C_{\alpha_{k}-1}=\binom{2n-1}{n-2}$$
which holds for every $n\in \mathbb{N}$ such that $n\geq 2$. Indeed, the left hand side of the above equation counts all matchings in $\mu_{n}(1212)$ by deleting the rightmost edge, then counting the resulting $1212-$avoiding matchings according to the number of factors. As an immediate consequence of Proposition $\ref{W}$,  $\ref{J}$ and   $\ref{1212}$ we deduce the following.

\begin{theorem} \label{T}
Let $\sigma\in \{1212,1221\}$ and let $\tau$ be a matching. Then, for $n\geq 2$,
$$|\mathcal{M}_{n}(\sigma(\tau+2))|=C_{n}+\sum_{\ell=2}^{n}\sum_{k=0}^{n-\ell}\binom{2\ell-1}{\ell-2}\binom{2\ell-1+k}{k}\binom{2(n-\ell)-k}{k}k!|\mathcal{M}_{n-\ell-k}(\tau)|.$$
\end{theorem}

Specializing $\tau$ in Theorem $\ref{T}$, we are able to enumerate a couple of new classes of matchings avoiding a single pattern (see also Figure $\ref{TB}$).

\begin{corollary}\label{E}
Let $n\in \mathbb{N}$, with $n\geq 2$, and $\sigma\in \{1212,1221\}$.
\begin{itemize}
\item[(i)] If $\tau\in \{1212,1221\}$, then  $$|\mathcal{M}_{n}(\sigma(\tau+2))|=C_{n}+\sum_{\ell=2}^{n}\sum_{k=0}^{n-\ell}\binom{2\ell-1}{\ell-2}\binom{2\ell+k-1}{k}\binom{2n-2\ell-k}{k}k!C_{n-\ell-k}.$$
\item[(ii)] If $\tau\in \{123123,123321\}$, then $$|\mathcal{M}_{n}(\sigma(\tau+2))|=C_{n}+\sum_{\ell=2}^{n}\sum_{k=0}^{n-\ell}\binom{2\ell-1}{\ell-2}\binom{2\ell+k-1}{k}\binom{2n-2\ell-k}{k}k!(C_{n-\ell-k}C_{n-\ell-k+2}-C_{n-\ell-k+1}^{2}).$$
\end{itemize}
\end{corollary}

\begin{figure}\label{TB}
\begin{center}
\begin{tabular}{|c|c|c|}
  \hline
  $n$ & $\mathcal{M}_{n}(12123434)$ &  $\mathcal{M}_{n}(1212345345)$ \\
  \hline
  1 & 1 & 1\\2 & 3 & 3\\ 3 & 15 & 15\\4 & 104 & 105\\5 & 910 & 944\\6 & 9503 & 10341\\7 & 114317 & 133132\\8 & 1547124 & 1961919\\ 9 & 23169162 & 32441303\\ 10 & 379308106 & 592718236 \\
  \hline
\end{tabular}
\end{center}
\caption{The first terms of the sequences of Corollary $\ref{E}$. This sequences are not recorded in $\cite{S}$.}
\end{figure}

\subsection{The lifting of a pattern}

 In this section we investigate classes of matchings avoiding the lifting of a given matching $\sigma$. The enumeration of such classes seems to be a hard problem in general, since a special instance of it is the enumeration of matchings avoiding the pattern $123231$, which is the lifting of $1212$, and it was remarked in Section $\ref{Intro}$ that this is likely to be a hard problem. However, if we impose additional constraints, namely the avoidance of a special pattern $\chi$ and its reversal $\overline{\chi}$, the description of the structure of matchings avoiding the lifting of $\sigma$ becomes more accessible. We start by fixing some preliminary definitions. Let $e$ and $f$ be any two edges of $\sigma$. We say that $e$ is $\emph{nested}$ in $f$ when $\min(f)<\min(e)$ and $\max(e)<\max(f)$. We say that $e$ is a $\emph{nested}$ $\emph{edge}$ when it is nested in some edge of $\sigma$ and that $e$ is a $\emph{top}$ $\emph{edge}$ otherwise. The pattern of $\sigma$ consisting of all the nested edges of $\sigma$ will be called the $\emph{core}$ of $\sigma$ and the pattern of $\sigma$ consisting of all the top edges of $\sigma$ will be called the $\emph{roof}$ of $\sigma$. Note that, by definition, the roof of $\sigma$ is a nonnesting matching. We say that a matching is $\emph{connected}$ when it is nonempty and it is not the juxtaposition of two nonempty matchings. Let $S$ be a set of matchings and $n\in \mathbb{N}$, we denote by $\mathcal{M}^{*}(S)$ the class of all connected matchings, by $\mathcal{M}_{n}^{*}(S)$ the set of matchings in $\mathcal{M}^{*}(S)$ of order $n$ and by $\mathcal{M}^{*}(S,z)$ the  generating function of $\mathcal{M}^{*}(S)$. In the following we will make some use of the so-called symbolic method, borrowing some standard constructions and notations from $\cite{FS}$, such as disjoint union, cartesian product and composition of combinatorial classes (in particular, the operator $\mathrm{Seq}$), which will allow us to easily translate combinatorial descriptions into  generating functions.

\bigskip

\begin{remark}\label{R}
 Note that, by definition, $\mathcal{M}(S)=\mathrm{Seq}(\mathcal{M}^{*}(S))$ and therefore $\mathcal{M}(S,z)=\frac{1}{1-\mathcal{M}^{*}(S,z)}$. In particular, $C(z)=\mathcal{M}(1221,z)=\frac{1}{1-\mathcal{M}^{*}(1221,z)}$, which leads to $$\mathcal{M}^{*}(1221,z)=\frac{C(z)-1}{C(z)}=\frac{zC(z)^{2}}{C(z)}=zC(z)$$ which means that, for every $n\in \mathbb{N}^{*}$, there are $C_{n-1}$ connected nonnesting matchings of order $n$.
\end{remark}

\bigskip

We are now in a position to state and prove the main result of this section.

\begin{theorem}\label{L} Let $\sigma$ be a connected matching and set $\chi=123132$, so that $\overline{\chi}=123213$. Then $$\mathcal{M}(1(\sigma+1)1,\chi,\overline{\chi},z)=\frac{1}{1-z\mathcal{M}(\sigma,\chi,\overline{\chi},z)C(z \mathcal{M}(\sigma,\chi,\overline{\chi},z)^{2})}.$$
\end{theorem}

\proof
Let $n\in \mathbb{N}^{*}$, $\lambda\in \mathcal{M}^{*}_{n}(1(\sigma+1)1,\chi,\overline{\chi})$ and $m$ be the order of its roof. The matching $\lambda$ is required to avoid both $\chi$ and $\overline{\chi}$, that are the matchings represented by the following linear chord diagrams:

\begin{center}
\scalebox{1.3}{\begin{tikzpicture}
\node at (0,0){
\scalebox{0.3}{
\begin{tikzpicture}
\draw [fill] (1,0) circle [radius=0.1];
\draw [fill] (2,0) circle [radius=0.1];
\draw [fill] (3,0) circle [radius=0.1];
\draw [fill] (4,0) circle [radius=0.1];
\draw [fill] (5,0) circle [radius=0.1];
\draw [fill] (6,0) circle [radius=0.1];
\draw[thick] (1,0) to [out=45,in=135] (4,0);
\draw[thick] (2,0) to [out=45,in=135] (6,0);
\draw[thick] (3,0) to [out=45,in=135] (5,0);
\end{tikzpicture}}};

\node at (4,0) {\scalebox{0.3}{
\begin{tikzpicture}
\draw [fill] (1,0) circle [radius=0.1];
\draw [fill] (2,0) circle [radius=0.1];
\draw [fill] (3,0) circle [radius=0.1];
\draw [fill] (4,0) circle [radius=0.1];
\draw [fill] (5,0) circle [radius=0.1];
\draw [fill] (6,0) circle [radius=0.1];
\draw[thick] (1,0) to [out=45,in=135] (5,0);
\draw[thick] (2,0) to [out=45,in=135] (4,0);
\draw[thick] (3,0) to [out=45,in=135] (6,0);
\end{tikzpicture}}};
\end{tikzpicture}}
\end{center}
This means that every nested edge of $\lambda$ is forced to never cross a top edge of $\lambda$. Therefore the core of $\lambda$ can be decomposed as the juxtaposition of $2m-1$ (possibly empty) matchings $\lambda_{1},...,\lambda_{2m-1}\in \mathcal{M}(\sigma,\chi,\overline{\chi})$, moreover  the occurrences of these factors in $\lambda$ are separated by the vertices of the top edges of $\lambda$. Conversely, every matching constructed as above belongs to the class $\mathcal{M}^{*}(1(\sigma+1)1,\chi,\overline{\chi})$, because $\sigma$ is connected and so no occurrence of $\sigma$ can show up by juxtaposing two patterns in the class $\mathcal{M}(\sigma,\chi,\overline{\chi})$. Thus   $\mathcal{M}^{*}_{n}(1(\sigma+1)1,\chi,\overline{\chi})$ is the set of matchings obtained by choosing some $m\in \mathbb{N}^{*}$ and a matching in $\mathcal{M}_{m}^{*}(1221)$, then replacing its edges other than the rightmost one with $(\{\scalebox{0.4}{
\begin{tikzpicture}
\draw [fill] (1,0) circle [radius=0.1];
\draw [fill] (2,0) circle [radius=0.1];
\draw[thick] (1,0) to [out=45+15,in=135-15] (2,0);
\end{tikzpicture}}\}\times \mathcal{M}^{*}(\sigma,\chi,\overline{\chi})^{2})-$structures and the rightmost edge with a $(\{\scalebox{0.4}{
\begin{tikzpicture}
\draw [fill] (1,0) circle [radius=0.1];
\draw [fill] (2,0) circle [radius=0.1];
\draw[thick] (1,0) to [out=45+15,in=135-15] (2,0);
\end{tikzpicture}}\}\times \mathcal{M}^{*}(\sigma,\chi,\overline{\chi}))-$structure. An instance of this decomposition is illustrated in the following figure when the roof is $121323$.

 \begin{center}
\scalebox{0.8}{
\begin{tikzpicture}
\draw[thick, green!50!black] (1,0) to [out=45,in=135] (7,0);
\draw[thick, green!50!black] (4,0) to [out=45,in=135] (13,0);
\draw[thick] (10,0) to [out=45,in=135] (16,0);

\node at (1,-0.5) {\scalebox{0.8}{\begin{tikzpicture}
   \node (pict) at (0,0) {
   \scalebox{1}{
\begin{tikzpicture}
\draw [fill, red] (1,0) circle [radius=0.1];
\draw [fill, red] (1.5,0) circle [radius=0.1];
\draw [fill, red] (3,0) circle [radius=0.1];
\draw[dashed,red, fill=red!25, pattern=crosshatch dots, pattern color=red] (1.5,0) to [out=45+40,in=135-40] (3,0) -- (1.5,0);
\end{tikzpicture}}
    };
   \node[draw,thick,fit=(pict),rounded corners=.55cm,]    {};
\end{tikzpicture}}};

\node at (4,-0.5) {\scalebox{0.8}{\begin{tikzpicture}
   \node (pict) at (0,0) {
   \scalebox{1}{
\begin{tikzpicture}
\draw [fill, red] (1,0) circle [radius=0.1];
\draw [fill, red] (1.5,0) circle [radius=0.1];
\draw [fill, red] (3,0) circle [radius=0.1];
\draw[dashed,red, fill=red!25, pattern=crosshatch dots, pattern color=red] (1.5,0) to [out=45+40,in=135-40] (3,0) -- (1.5,0);
\end{tikzpicture}}
    };
   \node[draw,thick,fit=(pict),rounded corners=.55cm,]    {};
\end{tikzpicture}}};

\node at (7,-0.5) {\scalebox{0.8}{\begin{tikzpicture}
   \node (pict) at (0,0) {
   \scalebox{1}{
\begin{tikzpicture}
\draw [fill, red] (1,0) circle [radius=0.1];
\draw [fill, red] (1.5,0) circle [radius=0.1];
\draw [fill, red] (3,0) circle [radius=0.1];
\draw[dashed,red, fill=red!25, pattern=crosshatch dots, pattern color=red] (1.5,0) to [out=45+40,in=135-40] (3,0) -- (1.5,0);
\end{tikzpicture}}
    };
   \node[draw,thick,fit=(pict),rounded corners=.55cm,]    {};
\end{tikzpicture}}};

\node at (10,-0.5) {\scalebox{0.8}{\begin{tikzpicture}
   \node (pict) at (0,0) {
   \scalebox{1}{
\begin{tikzpicture}
\draw [fill, red] (1,0) circle [radius=0.1];
\draw [fill, red] (1.5,0) circle [radius=0.1];
\draw [fill, red] (3,0) circle [radius=0.1];
\draw[dashed,red, fill=red!25, pattern=crosshatch dots, pattern color=red] (1.5,0) to [out=45+40,in=135-40] (3,0) -- (1.5,0);
\end{tikzpicture}}
    };
   \node[draw,thick,fit=(pict),rounded corners=.55cm,]    {};
\end{tikzpicture}}};

\node at (13,-0.5) {\scalebox{0.8}{\begin{tikzpicture}
   \node (pict) at (0,0) {
   \scalebox{1}{
\begin{tikzpicture}
\draw [fill, red] (1,0) circle [radius=0.1];
\draw [fill, red] (1.5,0) circle [radius=0.1];
\draw [fill, red] (3,0) circle [radius=0.1];
\draw[dashed,red, fill=red!25, pattern=crosshatch dots, pattern color=red] (1.5,0) to [out=45+40,in=135-40] (3,0) -- (1.5,0);
\end{tikzpicture}}
    };
   \node[draw,thick,fit=(pict),rounded corners=.55cm,]    {};
\end{tikzpicture}}};

\node at (16,-0.45) {\scalebox{1.3}{\begin{tikzpicture}
   \node (pict) at (0,0) {
   \scalebox{1}{
\begin{tikzpicture}
\draw [fill, red] (1,0) circle [radius=0.06];
\end{tikzpicture}}
    };
   \node[draw,line width=0.6,fit=(pict),rounded corners=0.30cm,]    {};
\end{tikzpicture}}};

\end{tikzpicture}}
\end{center}

It follows that the combinatorial class $\mathcal{M}^{*}(1(\sigma+1)1,\chi,\overline{\chi})$  is isomorphic to the combinatorial class $$\{\scalebox{0.4}{
\begin{tikzpicture}
\draw [fill] (1,0) circle [radius=0.1];
\draw [fill] (2,0) circle [radius=0.1];
\draw[thick] (1,0) to [out=45+15,in=135-15] (2,0);
\end{tikzpicture}}\}\times \mathcal{M}(\sigma,\chi,\overline{\chi})\times \sum_{m\geq 1}\mathcal{M}_{m}^{*}(1221)\times (\{\scalebox{0.4}{
\begin{tikzpicture}
\draw [fill] (1,0) circle [radius=0.1];
\draw [fill] (2,0) circle [radius=0.1];
\draw[thick] (1,0) to [out=45+15,in=135-15] (2,0);
\end{tikzpicture}}\}\times \mathcal{M}(\sigma,\chi,\overline{\chi})^{2})^{m-1}$$
 and this isomorphism immediately translates into the following expression for the generating function
\begin{align*}
\mathcal{M}^{*}(1(\sigma+1)1,\chi,\overline{\chi},z) & = z\mathcal{M}(\sigma,\chi,\overline{\chi},z)\sum_{m\geq 1}[z^{m}](zC(z)) (z\mathcal{M}(\sigma,\chi,\overline{\chi},z)^{2})^{m-1}\\
 & =z\mathcal{M}(\sigma,\chi,\overline{\chi},z)\sum_{m\geq 0}C_{m} (z\mathcal{M}(\sigma,\chi,\overline{\chi},z)^{2})^{m}\\
 & =z\mathcal{M}(\sigma,\chi,\overline{\chi},z)C(z\mathcal{M}(\sigma,\chi,\overline{\chi},z)^{2}).\\
\end{align*}
Now the claim follows from the above Remark.
\endproof

Note that, at least in principle,  iterating Theorem $\ref{L}$ allows us to find expressions for the generating function of $\mathcal{M}(12...k(\sigma+k)k...21,\chi,\overline{\chi})$ in terms of the generating function of $\mathcal{M}(\sigma,\chi,\overline{\chi})$, for every $k\in \mathbb{N}^{*}$. As an immediate application, we are able to compute the generating function of two classes of matchings avoiding three patterns of order three.

\begin{corollary}\label{123231}
The following equality holds
 $$\mathcal{M}(123231,123132,123213,z)=\mathcal{M}(123321,123132,123213,z)=\frac{1}{1-zC(z)C(C(z)-1)}$$ and $|\mathcal{M}_{n}(123231,123132,123213)|=|\mathcal{M}_{n}(123321,123132,123213)|$ is the $n^{th}$ term of sequence $\mathrm{A125188}$ in $\cite{S}$.
\end{corollary}

\proof Let $\sigma\in \{1212,1221\}$, then it follows from Theorem $\ref{L}$ that
$$\mathcal{M}(1(\sigma+1)1,\chi,\overline{\chi},z)=\frac{1}{1-z\mathcal{M}(\sigma,\chi,\overline{\chi},z)C(\mathcal{M}(\sigma,\chi,\overline{\chi},z))}$$
where $\chi=123132$ and $\overline{\chi}=123213$. Moreover, $\mathcal{M}(\sigma,\chi,\overline{\chi},z)=\mathcal{M}(\sigma,z)=C(z)$ and the first claim follows. The generating function for sequence $\mathrm{A125188}$ can be found in $\cite{S}$ and is given by
$$f(z)=\frac{1+zC(z)-\sqrt{1-zC(z)-5z}}{2z(1+C(z))}$$
Applying the change of variable $y=zC(z)$, so that $z=y(1-y)$ and $C(z)=\frac{1}{1-y}$, some routine computations show that $f(z)=\frac{1}{1-zC(z)C(C(z)-1)}$, hence the second claim   also follows.
\endproof

Sequence $\mathrm{A125188}$ counts Dumont permutations of the first kind avoiding the patterns $2413$ and $4132$, but we have not been able to find any bijection with our classes of pattern avoiding matchings. Note that, for $\sigma\in \{1212,1221\}$, iterating Theorem $\ref{L}$ allows to prove that $\mathcal{M}(12...k(\sigma+k)k...21,\chi,\overline{\chi},z)$ is an algebraic function of $C(z)$, hence it is itself algebraic, for every $k\in \mathbb{N}^{*}$.

\section{Unlabeled pattern avoidance}\label{UPA}

In this section we introduce the notion of unlabeled matching, which provides a way to collect patterns  that are combinatorially equivalent, in a sense that is specified below. Given $n\in \mathbb{N}^{*}$, let $\gamma_{n}$ denote the $2n-$cycle $(1\ 2\ 3\ ...\ 2n)$ on $[2n]$ and let $\sigma$ and $\tau$ be two matchings of order $n$. We say that $\sigma$ and $\tau$ are $\emph{ciclically}$ $\emph{equivalent}$ when there exists $k\in [2n]$ such that $\{i,j\}\in \sigma$ if and only if $\{\gamma_{n}^{k}(i),\gamma_{n}^{k}(j)\}\in \tau$, for every $i,j\in [2n]$. In other words, two matchings are ciclically  equivalent when they have the same unlabeled circular chord diagram. An equivalence class of matchings is called an $\emph{unlabeled}$ $\emph{matching}$. For instance, $[112323]=\{112323,123231,123312,121233,121332,122313\}$. Thus, an unlabeled matching can be represented by an unlabeled circular chord diagram; for instance, the unlabeled matching $[112323]$ can be represented by the following unlabeled chord diagram

\begin{center}
\scalebox{0.25}{\begin{tikzpicture}[x=1.00mm, y=1.00mm, inner xsep=0pt, inner ysep=0pt, outer xsep=0pt, outer ysep=0pt]
\path[line width=0mm] (67.21,8.97) rectangle +(104.91,102.64);
\definecolor{L}{rgb}{0,0,0}
\path[line width=0.30mm, draw=L] (119.78,60.29) circle (49.32mm);
\path[line width=0.30mm, draw=L] (71.36,70.16) arc (-147:-147:69.96mm);
\path[line width=0.60mm, draw=L] (70.21,60.24) arc (-83:-23:80.19mm);
\definecolor{F}{rgb}{0,0,0}
\path[line width=0.60mm, draw=L, fill=F] (133.61,107.51) circle (2.50mm);
\path[line width=0.60mm, draw=L, fill=F] (70.21,60.01) circle (2.50mm);
\path[line width=0.60mm, draw=L] (105.25,107.65) arc (143:143:79.41mm);
\path[line width=0.60mm, draw=L] (105.72,107.42) arc (-157:-98:81.38mm);
\path[line width=0.60mm, draw=L, fill=F] (105.72,107.65) circle (2.50mm);
\path[line width=0.60mm, draw=L, fill=F] (169.12,59.23) circle (2.50mm);
\path[line width=0.60mm, draw=L] (157.36,31.25) arc (58:123:71.73mm);
\path[line width=0.60mm, draw=L, fill=F] (158.75,30.32) circle (2.50mm);
\path[line width=0.60mm, draw=L, fill=F] (80.81,30.55) circle (2.50mm);
\end{tikzpicture}}
\end{center}
Note that a matching avoids an unlabeled pattern if and only if its circular chord diagram avoids the unlabeled chord diagram of the pattern.

The unlabeled matchings of order $2$ are exactly $[1122]=\{1122,1221\}$ and $[1212]=\{1212\}$.
Note that, for every $n\in \mathbb{N}^{*}$, a matching  $\lambda$ of order $n$ avoids $[1122]$ if and only if it is permutational and nonnesting, hence $\mathcal{M}_{n}([1122])=\{123...n123...n\}$.
We thus have $\mathcal{M}([1212],z)=C(z)$ and $\mathcal{M}([1122],z)=\frac{1}{1-z}$. The unlabeled matchings of order $3$ are exactly five, namely:
\begin{itemize}
\item[]$[112323]=\{ 112323,123231,123312,121233,121332,122313 \},$
\item[]$[123132]=\{ 121323,123213,121323\},$
\item[]$[123321]=\{ 123321,122133,112332\},$
\item[]$[112233]=\{ 112233,122331\},$
\item[]$[123123]=\{ 123123\}.$
\end{itemize}
Clearly $|\mathcal{M}_{n}([123123])|=C_{n}C_{n+2}-C_{n+1}^{2}$. In this section we will work out explicit formulas to enumerate  $\mathcal{M}([112323])$ and  $\mathcal{M}([123132])$.

\begin{proposition} The generating function of matchings avoiding the unlabeled pattern $[112323]$ is given by  $$\mathcal{M}([112323],z)=C(z)+\frac{z^{2}}{(1-z)^{2}(1-2z)}$$ As a consequence, its coefficients have the following closed form: $$|\mathcal{M}_{n}([112323])|=C_{n}+2^{n}-n-1,$$ for $n\geq 2$.
\end{proposition}

\proof Clearly the noncrossing matchings in $\mathcal{M}_{n}([112323])$ are counted by the Catalan number $C_{n}$, hence it remains to count the crossing matchings in $\mathcal{M}_{n}([112323])$. Let  $\lambda$ be a crossing matching in $\mathcal{M}_{n}([112323])$. Let $\sigma$ denote the pattern of $\lambda$ consisting of all the edges intersecting the leftmost edge of $\lambda$ and let $\tau$ denote the pattern of $\lambda$ consisting of all the remaining edges. Note that $\sigma$ is nonempty, otherwise, since $\lambda$ is assumed to be crossing, there would be a pair of crossing edges that do not cross the leftmost edge of $\lambda$, thus forming an occurrence of $[112323]$. Assume that $\sigma$ contains $k$ edges, where $k\in [n-1]$.   Observe that $\sigma$ has to be permutational, because an occurrence of $1122$ in $\sigma$ should have at least one edge which does not cross the leftmost edge of $\lambda$, against the definition of $\sigma$. 
Moreover, since $\sigma$ avoids $[112323]$, the corresponding permutation has to avoid both the permutation patterns $231$ and $312$, therefore there are $|\mathcal{S}_{k}(231,312)|=2^{k-1}$ possible choices for $\sigma$.  Furthermore, $\tau$ must be noncrossing, otherwise any pair of crossing edges of $\tau$ together with the leftmost edge of $\lambda$ would form an occurrence of $[112323]$. Finally, using a similar argument, we deduce that each edge  of $\tau$ has to cross all the edges of $\sigma$. 
We can thus conclude that $\tau$ is the juxtaposition of two totally crossing matchings of order $n-k$ such that the leftmost one is nonempty. Hence there are exactly $n-k$ possibile choices for $\tau$. In other words, $\lambda$ has the form illustrated by the following linear chord diagram

\begin{center}
\scalebox{0.5}{\begin{tikzpicture}
\draw [fill] (1,0) circle [radius=0.1];
\draw [fill] (2,0) circle [radius=0.1];
\draw [fill] (4,0) circle [radius=0.1];
\draw [fill] (5,0) circle [radius=0.1];
\draw [fill] (7,0) circle [radius=0.1];
\draw [fill] (8,0) circle [radius=0.1];
\draw [fill] (10,0) circle [radius=0.1];
\draw [fill] (11,0) circle [radius=0.1];
\draw [fill] (14,0) circle [radius=0.1];
\draw [fill] (15,0) circle [radius=0.1];
\draw [fill] (17,0) circle [radius=0.1];
\draw [fill] (18,0) circle [radius=0.1];
\draw [fill] (20,0) circle [radius=0.1];
\draw [fill] (21,0) circle [radius=0.1];
\draw [fill] (23,0) circle [radius=0.1];
\draw [fill] (24,0) circle [radius=0.1];
\draw  (1,0) to [out=45,in=135] (11,0);
\draw  (2,0) to [out=45,in=135] (10,0);
\draw  (4,0) to [out=45,in=135] (8,0);
\draw  (14,0) to [out=45,in=135] (24,0);
\draw  (15,0) to [out=45,in=135] (23,0);
\draw  (17,0) to [out=45,in=135] (21,0);
\draw[dashed, pattern=crosshatch dots] (5,0) to [out=45,in=135] (20,0) -- (18,0) to [out=135,in=45] (7,0) -- (5,0) ;
\node at (3,0) {\scalebox{1}{$...$}};
\node at (9,0) {\scalebox{1}{$...$}};
\node at (16,0) {\scalebox{1}{$...$}};
\node at (22,0) {\scalebox{1}{$...$}};
\node at (12.5,2.6) {\scalebox{1.5}{$\sigma$}};
\node at (3.5,0.5) {\scalebox{1.5}{$\tau$}};
\end{tikzpicture}}
\end{center}

 From this characterization of $\mathcal{M}_{n}([112323])$, it follows that
$$|\mathcal{M}_{n}([112323])|=C_{n}+\sum_{k=1}^{n-1}(n-k)2^{k-1},$$
hence
\begin{align*}
\mathcal{M}([112323],z) &=C(z)+\sum_{n\geq 2}\sum_{k=1}^{n-1}(n-k)2^{k-1}z^{n} \\
 &=C(z)+z^{2}\sum_{n\geq 0}\sum_{k=0}^{n}(n-k+1) 2^{k}z^{n}\\
 &=C(z)+z^{2}\left(\sum_{n\geq 0}(n+1)z^{n}\right)\left(\sum_{n\geq 0}2^{n}z^{n}\right)\\
 &=C(z)+\frac{z^{2}}{(1-z)^{2}(1-2z)}.\\
\end{align*}
 Finally we can compute the partial fraction decomposition of $\frac{z^{2}}{(1-z)^{2}(1-2z)}$ to find
\begin{align*}
\frac{z^{2}}{(1-z)^{2}(1-2z)}& =z^{2}\left[\frac{-3+2z}{(1-z)^{2}}+\frac{4}{1-2z}\right]\\
 & =z^{2}\left[-3\sum_{n\geq 0}(n+1)z^{n}+2z\sum_{n\geq 0}(n+1)z^{n}+4\sum_{n\geq 0}2^{n}z^{n}\right]\\
 & =\sum_{n\geq 2}(2^{n}-n-1)z^{n},
\end{align*}
which proves the claim.

\endproof

The sequence enumerating $\mathcal{M}([112323])$ begins $1,1,1,3,9,25,68,189,...$ and it is not recorded in $\cite{S}$, however it is worth noting that $2^{n}-n-1$ is the $n^{th}$ Eulerian number (sequence $\mathrm{A000295}$ in $\cite{S}$).

Our last result concerns the unlabeled pattern $[123132]$, which is represented by the following unlabeled chord diagram:

\begin{center}
\scalebox{0.3}{
\begin{tikzpicture}[x=1.00mm, y=1.00mm, inner xsep=0pt, inner ysep=0pt, outer xsep=0pt, outer ysep=0pt]
\path[line width=0mm] (43.31,17.06) rectangle +(118.59,86.01);
\definecolor{L}{rgb}{0,0,0}
\path[line width=0.30mm, draw=L] (119.78,60.06) circle (40.12mm);
\path[line width=0.60mm, draw=L] (45.31,83.14) arc (-73:-73:36.17mm);
\path[line width=0.60mm, draw=L] (87.04,83.20) arc (-116:-64:74.45mm);
\path[line width=0.60mm, draw=L] (152.52,36.91) arc (64:116:74.35mm);
\path[line width=0.60mm, draw=L] (120.01,100.07) -- (120.01,20.06);
\definecolor{F}{rgb}{0,0,0}
\path[line width=0.30mm, draw=L, fill=F] (120.01,100.07) circle (2.00mm);
\path[line width=0.30mm, draw=L, fill=F] (87.04,83.47) circle (2.00mm);
\path[line width=0.30mm, draw=L, fill=F] (152.52,83.47) circle (2.00mm);
\path[line width=0.30mm, draw=L, fill=F] (86.81,36.90) circle (2.00mm);
\path[line width=0.30mm, draw=L, fill=F] (152.29,36.90) circle (2.00mm);
\path[line width=0.30mm, draw=L, fill=F] (120.01,20.06) circle (2.00mm);
\end{tikzpicture}}
\end{center}

It turns out that matchings avoiding $[123132]$ have a ternary tree structure and the following discussion is in fact devoted to describe a bijection between this class of matchings and ternary trees.  To this purpose, recall that, for $k\in \mathbb{N}^{*}$, a $\emph{k-ary}$ $\emph{tree}$ is an ordered rooted tree such that every node has at most $k$ children. Let $\mathcal{T}_{k}$ denote the combinatorial class of $k-$ary trees. Note that every $k-$ary tree is either empty or it can be decomposed as in the following figure:
\begin{center}
\begin{tikzpicture}
\draw [fill] (3,1) circle [radius=0.1];
\node at (1,-0.5) {$T_{1}$};
\node at (2,-0.5) {$T_{2}$};
\node at (3,-0.5) {$\hdots$};
\node at (4,-0.5) {$T_{k-1}$};
\node at (5,-0.5) {$T_{k}$};
\draw[thick] (1,-0.5) circle [radius=0.5];
\draw[thick] (2,-0.5) circle [radius=0.5];
\draw[thick] (4,-0.5) circle [radius=0.5];
\draw[thick] (5,-0.5) circle [radius=0.5];
\draw[thick] (3,1)--(1,0);
\draw[thick] (3,1)--(2,0);
\draw[thick] (3,1)--(4,0);
\draw[thick] (3,1)--(5,0);
\end{tikzpicture}
\end{center}
where $\bullet$ is the root and $T_{1},...,T_{k}\in \mathcal{T}_{k}$. Therefore the combinatorial classes $\mathcal{T}_{k}$ and $\{\emptyset\}+\{\bullet\}\times (\mathcal{T}_{k})^{k}$ are isomorphic and the isomorphism  translates into the functional equation  $\mathcal{T}_{k}(z)=1+z\mathcal{T}_{k}(z)^{k}$ for the generating function $\mathcal{T}_{k}(z)$ of the class $\mathcal{T}_{k}$. This equation can be classically solved by Lagrange's inversion as follows
$$[z^{n}]\mathcal{T}_{k}(z)=[z^{n}](\mathcal{T}_{k}(z)-1)=\frac{1}{n}[w^{n-1}](1+w)^{kn}=$$ $$\frac{1}{n}[w^{n-1}]\sum_{i=0}^{kn}\binom{kn}{i}w^{i}=\frac{1}{n}\binom{kn}{n-1}=\frac{1}{(k-1)n+1}\binom{kn}{n}$$
In particular, when $k=3$, we thus get $$[z^{n}]\mathcal{T}_{3}(z)=\frac{1}{2n+1}\binom{3n}{n}$$
for every $n\in \mathbb{N}^{*}$.

Now we recursively define a map $\varphi:\{\emptyset\}+\{\bullet\}\times(\mathcal{T}_{3})^{3}\longrightarrow \mathcal{M}([123132])$ as follows. Set $\varphi(\emptyset)=\emptyset$; furthermore, for every $(T_{1},T_{2},T_{3})\in (\mathcal{T}_{3})^{3}$, let $\varphi(\bullet,T_{1},T_{2},T_{3})$ be the matching whose linear chord diagram $\Gamma$  is constructed as follows:

\begin{itemize}
\item[1.] Denote with $\Gamma^{(i)}$ the linear chord diagram of $\varphi(T_{i})$, for every $i\in \{1,2,3\}$.
\item[2.] If $\Gamma^{(1)}$ is empty, then:
\begin{itemize}
\item[1.] draw a vertex $\ell'$ on the left of a vertex $r'$ and connect them with an edge;
\item[2.] draw $\Gamma^{(2)}$ between $\ell'$ and $r'$ and $\Gamma^{(3)}$ to the right of $r'$.
\end{itemize}
\item[3.] If $\Gamma^{(1)}$ is nonempty, then:
\begin{itemize}
\item[1.] let $\ell$ and $r$ denote the left and right vertex of the leftmost edge of $\Gamma^{(1)}$, respectively; draw two vertices $\ell'$ and $r'$ to the left of $\ell$ and $r$, respectively, and connect them with an edge.
\item[2.] draw  $\Gamma_{2}$ between $\ell'$ and $\ell$ and $\Gamma_{3}$ between $r'$ and $r$.
\end{itemize}
\end{itemize}

 In other words, the map $\varphi$ can be represented by the following diagram:

 \begin{center}
\scalebox{1.7}{\includegraphics[width=0.5\textwidth]{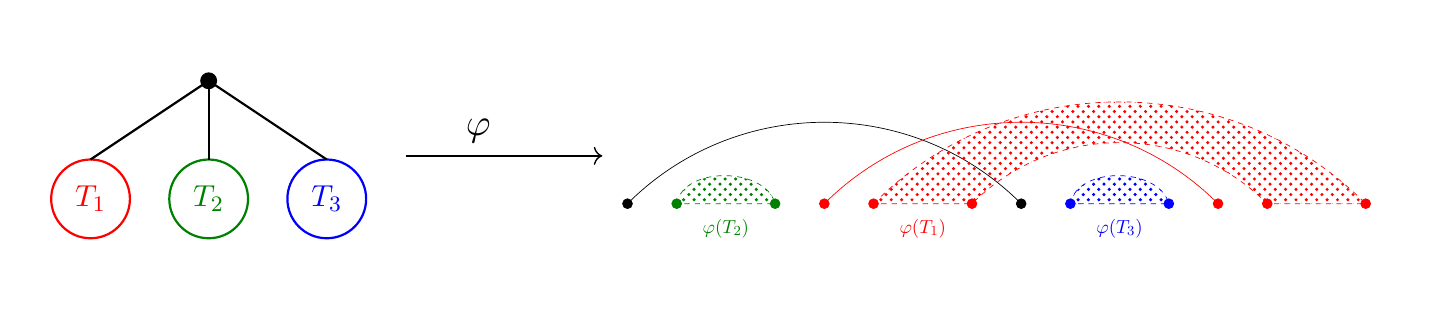}}
 \end{center}

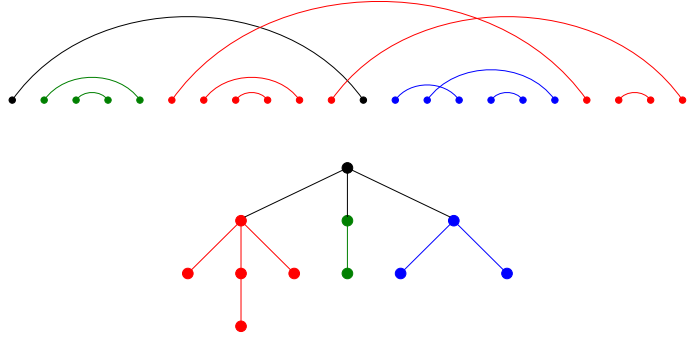
\begin{figure}
\begin{center}\scalebox{0.7}{
\begin{tikzpicture}
\node at (0,0)
{\begin{tikzpicture}
\draw[fill] (3,3) circle (1mm);

\draw[fill, red] (1,2) circle (1mm);
\draw[fill,green!50!black] (3,2) circle (1mm);
\draw[fill, blue] (5,2) circle (1mm);

\draw[fill,red] (0,1) circle (1mm);
\draw[fill,red] (1,1) circle (1mm);
\draw[fill,red] (2,1) circle (1mm);
\draw[fill,green!50!black] (3,1) circle (1mm);
\draw[fill, blue] (4,1) circle (1mm);
\draw[fill, blue] (6,1) circle (1mm);

\draw[fill,red] (1,0) circle (1mm);

\draw (3,3)--(1+0.07,2+0.07);
\draw[red] (1,2)--(0,1);
\draw[red] (1,2)--(1,1)--(1,0);
\draw[red] (1,2)--(2,1);
\draw (3,3)--(3,2+0.07);
\draw[green!50!black] (3,2)--(3,1);
\draw (3,3)--(5-0.07,2+0.07);
\draw[blue] (5,2)--(4,1);
\draw[blue] (5,2)--(6,1);

\end{tikzpicture}};

\node at (0,4) {
\scalebox{0.6}{\begin{tikzpicture}
\draw[fill] (1,0) circle (1mm);
\draw[fill, green!50!black] (2,0) circle (1mm);
\draw[fill, green!50!black] (3,0) circle (1mm);
\draw[fill, green!50!black] (4,0) circle (1mm);
\draw[fill, green!50!black] (5,0) circle (1mm);
\draw[fill, red] (6,0) circle (1mm);
\draw[fill, red] (7,0) circle (1mm);
\draw[fill, red] (8,0) circle (1mm);
\draw[fill, red] (9,0) circle (1mm);
\draw[fill, red] (10,0) circle (1mm);
\draw[fill, red] (11,0) circle (1mm);
\draw[fill] (12,0) circle (1mm);
\draw[fill, blue] (13,0) circle (1mm);
\draw[fill, blue] (14,0) circle (1mm);
\draw[fill, blue] (15,0) circle (1mm);
\draw[fill, blue] (16,0) circle (1mm);
\draw[fill, blue] (17,0) circle (1mm);
\draw[fill, blue] (18,0) circle (1mm);
\draw[fill, red] (19,0) circle (1mm);
\draw[fill, red] (20,0) circle (1mm);
\draw[fill, red] (21,0) circle (1mm);
\draw[fill, red] (22,0) circle (1mm);

\draw (1,0) to[out=45+10 ,in=135-10](12,0);
\draw[green!50!black ] (2,0) to[out=45+10 ,in=135-10](5,0);
\draw[green!50!black ] (3,0) to[out=45+10 ,in=135-10](4,0);
\draw[red] (6,0) to[out=45+10 ,in=135-10](19,0);
\draw[red ] (7,0) to[out=45+10 ,in=135-10](10,0);
\draw[red ] (8,0) to[out=45+10 ,in=135-10](9,0);
\draw[red ] (11,0) to[out=45+10 ,in=135-10](22,0);
\draw[blue ] (13,0) to[out=45+10 ,in=135-10](15,0);
\draw[blue ] (14,0) to[out=45+10 ,in=135-10](18,0);
\draw[blue ] (16,0) to[out=45+10 ,in=135-10](17,0);
\draw[red ] (20,0) to[out=45+10 ,in=135-10](21,0);
\end{tikzpicture}}};
\end{tikzpicture}}
\end{center}
\caption{The linear chord diagram of a matching with semilength 11 avoiding the unlabeled pattern $[123132]$ and the corresponding $3-$ary tree.}
\end{figure}

Conversely, define recursively a map $\psi:\mathcal{M}([123132])\longrightarrow \mathcal{T}_{3}$  as follows. Set $\psi(\emptyset)=\emptyset$. For every $\lambda\in \mathcal{M}([123132])\setminus\{\emptyset\}$, let $\psi(\lambda)$ be the ternary tree defined as follows:
\begin{itemize}
\item[1.] Suppose that the leftmost edge of $\lambda$ does not cross any other edge. In this case, denote by $\lambda_{2}$ the pattern of $\lambda$ consisting of all the edges of $\lambda$ which are nested below the leftmost edge of $\lambda$ and denote by $\lambda_{3}$ the pattern of $\lambda$ consisting of all the remaining edges of $\lambda$ other than the leftmost edge. We then define $\psi(\lambda)=(\bullet,\emptyset,\psi(\lambda_{2}),\psi(\lambda_{3}))$.
\item[2.] Suppose that the leftmost edge $\ell$ of $\lambda$ crosses some other edge of $\lambda$ and let $\ell'$ denote the leftmost edge of $\lambda$ among those crossed by $\ell$. Let $\lambda_{2}$ denote the pattern of $\lambda$ consisting of all $e\in \lambda$ such that $\min(\ell)<\min(e)<\max(e)<\min(\ell')$ and let $\lambda_{3}$ denote the pattern of $\lambda$ consisting of all $e\in \lambda$ such that $\max(\ell)<\min(e)<\max(e)<\max(\ell')$. Finally, let $\lambda_{1}$ denote the pattern of $\lambda$ consisting of all the remaining edges of $\lambda$ other than $\ell$. We then define $\psi(\lambda)=(\bullet,\psi(\lambda_{1}),\psi(\lambda_{2}),\psi(\lambda_{3}))$.  
\end{itemize}

\begin{proposition} The maps $\varphi$ and $\psi$ are well defined mutually inverse bijections.  In particular
$$|\mathcal{M}_{n}([123132])|=\frac{1}{2n+1}\binom{3n}{n},$$
for every $n\in \mathbb{N}$.
\end{proposition}

\proof The main thing we have to prove is that $\varphi$ is well defined. Denote by $\Lambda$ the unlabeled chord diagram of $[123132]$. Given $T=(\bullet,T_{1},T_{2},T_{3})\in \{\bullet\}\times (\mathcal{T}_{3})^{3}$, we now prove (by induction hypothesis on the number of nodes of $T$) that $\varphi(T)\in \mathcal{M}([123132])$. Using the same notation as in the definition of $\varphi$, we first observe that (by induction hypothesis) there is no occurrence of $\Lambda$ in $\Gamma^{(1)}$. Furthermore, the edge $\{\ell',r'\}$ cannot be involved in any occurrence of $\Lambda$. Suppose in fact that $\Lambda_{0}$ is an occurrence of $\Lambda$ involving $\{\ell',r'\}$. If $\Gamma^{(1)}$ is nonempty, then it is not difficult to realize that the leftmost edge $\{\ell,r\}$ of $\Gamma^{(1)}$ cannot occur in $\Lambda_{0}$ (this is due to the choice of the specific  pattern $\Lambda$). Thus we can replace $\{\ell',r'\}$ with $\{\ell,r\}$ in $\Lambda_{0}$ to get an occurrence of $\Lambda$ in $\Gamma^{(1)}$, which is a contradiction. On the other hand, if $\Gamma^{(1)}$ is empty, it is easy to check that $\{\ell',r'\}$ cannot belong to any occurrence of $\Lambda$ in $\varphi(T)$. Finally, no edge in $\Gamma^{(2)}$ or $\Gamma^{(3)}$ can be involved in an occurrence of $\Lambda$, because both $\varphi(T_{2})$ and $\varphi(T_{3})$ avoid $[123132]$ (by induction) and each of the edges of their chord diagrams does not cross any of the remaining edges of $\varphi(T)$. To conclude, it suffices to prove that $\varphi$ and $\psi$ are mutually inverse, which is immediate by their construction.
\endproof

\section{Combinatorics of intervals: preliminary results}\label{Int}

Another important topic that deserves to be investigated is the combinatorics of the intervals of the matching pattern poset.
In this sense, typical questions concern the counting of elements or, more generally, the enumeration of (saturated) chains of a given interval.
Another important line of research is the computation of the M\"obius function.
These are problems that have been classically studied for many combinatorial posets, such as Bruhat orders \cite{T} and Tamari lattices \cite{CCP,F}.
In this section we just scratch the surface of this vast subject,
by proposing a couple of relatively simple results concerning the enumeration of intervals of the form
$[\begin{tikzpicture}[scale=0.25] \draw (0,0) [fill] circle (.1); \draw (1,0) [fill] circle (.1); \draw[thick] (1,0) arc (0:180:0.5); \end{tikzpicture}\, ,\tau ]$, when $\tau$ has a specific form. In particular, in all the cases we will consider $\tau$ will be noncrossing.

\bigskip

Given a matching $\tau$, we say that an edge of $\tau$ is \emph{small} whenever its vertices are consecutive integers.
If $\tau (n,k)$ is a noncrossing matching of size $n$ having $k$ small edges, what is the cardinality of the interval
$[\begin{tikzpicture}[scale=0.25] \draw (0,0) [fill] circle (.1); \draw (1,0) [fill] circle (.1); \draw[thick] (1,0) arc (0:180:0.5); \end{tikzpicture}\, ,\tau (n,k)]$?
This may be a difficult problem in general. Here we address only a few very simple cases.

First of all, it is immediate to see that:
\begin{itemize}
\item $|[\begin{tikzpicture}[scale=0.25] \draw (0,0) [fill] circle (.1); \draw (1,0) [fill] circle (.1); \draw[thick] (1,0) arc (0:180:0.5); \end{tikzpicture}\, ,\tau (n,0)]|=0$,
    for all $\tau (n,0)$ (since there are no noncrossing matchings having no small edges);
\item $|[\begin{tikzpicture}[scale=0.25] \draw (0,0) [fill] circle (.1); \draw (1,0) [fill] circle (.1); \draw[thick] (1,0) arc (0:180:0.5); \end{tikzpicture}\, ,\tau (n,1)]|=n$,
    for all $\tau (n,1)$ (since, in this case, the interval is a chain having $n$ elements, which are all totally nesting matchings).
\end{itemize}

When $k=2$, the generic noncrossing matching having 2 small edges has the following form:

\bigskip

\begin{tikzpicture}[scale=0.5]
\draw (0,0) [fill] circle (.1); \draw (1,0) [fill] circle (.1); \draw (2,0) [fill] circle (.1); \draw (3,0) [fill] circle (.1);
\draw (4,0) [fill] circle (.1); \draw (5,0) [fill] circle (.1);
\draw[thick] (2,0) arc (0:180:0.5); \draw[thick] (4,0) arc (0:180:0.5); \draw[thick] (5,0) arc (0:180:2.5);
\node[above] at (1.5,0.5) {$r$}; \node[above] at (3.5,0.5) {$s$}; \node[above] at (2.5,2.5) {$k$};
\end{tikzpicture}

\bigskip

\noindent where an edge labeled $x$ stands for a totally nesting matching having $x$ edges.
In words, the above matching is the juxtaposition of two totally nesting matchings having $r$ and $s$ edges, respectively,
enclosed in a totally nesting matching having $k$ edges.
In order to have easy inline notations, such a matching will be denoted $\mathbf{k}(\mathbf{r};\mathbf{s})$.
Assuming w.l.o.g. that $r\geq s$, it is easy to see that
$[\begin{tikzpicture}[scale=0.25] \draw (0,0) [fill] circle (.1); \draw (1,0) [fill] circle (.1); \draw[thick] (1,0) arc (0:180:0.5); \end{tikzpicture}\, ,\mathbf{k}(\mathbf{r};\mathbf{s})]$
contains $r+k$ matchings having 1 small edge and $rs(k+1)$ matchings having 2 small edges. Therefore
$|[\begin{tikzpicture}[scale=0.25] \draw (0,0) [fill] circle (.1); \draw (1,0) [fill] circle (.1); \draw[thick] (1,0) arc (0:180:0.5); \end{tikzpicture}\, ,\mathbf{k}(\mathbf{r};\mathbf{s})]|=r+k+rs(k+1)$.

When $k=3$, again w.l.o.g., the generic matching $\tau (n,3)$ has the form

\bigskip

\begin{tikzpicture}[scale=0.5]
\draw (0,0) [fill] circle (.1); \draw (1,0) [fill] circle (.1); \draw (2,0) [fill] circle (.1); \draw (3,0) [fill] circle (.1);
\draw (4,0) [fill] circle (.1); \draw (5,0) [fill] circle (.1); \draw (6,0) [fill] circle (.1); \draw (7,0) [fill] circle (.1);
\draw (8,0) [fill] circle (.1); \draw (9,0) [fill] circle (.1);
\draw[thick] (3,0) arc (0:180:0.5); \draw[thick] (5,0) arc (0:180:0.5); \draw[thick] (8,0) arc (0:180:0.5);
\draw[thick] (6,0) arc (0:180:2.5); \draw[thick] (9,0) arc (0:180:4.5);
\node[above] at (2.5,0.5) {$a$}; \node[above] at (4.5,0.5) {$b$}; \node[above] at (7.5,0.5) {$c$};
\node[above] at (3.5,2.5) {$h$}; \node[above] at (4.5,4.5) {$k$};
\end{tikzpicture}

\bigskip

Similarly as before, we denote the above matching with $\mathbf{k}(\mathbf{h}(\mathbf{a};\mathbf{b});\mathbf{c})$.
We can count the elements of $[\begin{tikzpicture}[scale=0.25] \draw (0,0) [fill] circle (.1); \draw (1,0) [fill] circle (.1); \draw[thick] (1,0) arc (0:180:0.5); \end{tikzpicture}\, ,\mathbf{k}(\mathbf{h}(\mathbf{a};\mathbf{b});\mathbf{c})]$ with respect to the number of small edges.
\begin{itemize}
\item In order to count the number $\chi_1$ of matchings having one small edge,
we have to understand how many edges the largest totally nesting matching smaller than $\mathbf{k}(\mathbf{h}(\mathbf{a};\mathbf{b});\mathbf{c})$ has. To construct such a matching, we take the $k$ external edges, and add the largest number between $c$ and $h+\max (a,b)$. Therefore $\chi_1 =\max (k+h+a,k+h+b,k+c)$.
\item Matching having two edges can be obtained in two different ways from $\mathbf{k}(\mathbf{h}(\mathbf{a};\mathbf{b});\mathbf{c})$.
First, we can remove the totally nesting matchings having $c$, thus obtaining the matching $(\mathbf{h+k})(\mathbf{a};\mathbf{b})$,
which has $ab(h+k+1)$ matchings with 2 small edges below.
The second option is to remove one of the two totally nesting matchings with $a$ and $b$ edges, and precisely the smaller one,
thus obtaining the matching $\mathbf{k}((\mathbf{\max (a,b)+h});\mathbf{c})$, which has $(\max (a,b)+h)c(k+1)$ matchings with 2 small edges below.
However, there are matchings in common in the two above cases, which causes an overcount.
Indeed, the matchings which can be obtained in both the above cases are precisely those lying below $\mathbf{k}(\mathbf{a};\mathbf{\min (b,c)})$
and having 2 small edges,
which are $a\cdot \min (b,c)\cdot (k+1)$.
From the above consideration, we can write the total number $\chi_2$ of elements of the interval
$[\begin{tikzpicture}[scale=0.25] \draw (0,0) [fill] circle (.1); \draw (1,0) [fill] circle (.1); \draw[thick] (1,0) arc (0:180:0.5); \end{tikzpicture}\, ,\mathbf{k}(\mathbf{h}(\mathbf{a};\mathbf{b});\mathbf{c})]$ having 2 small edges, which is
$\chi_2 =ab(h+k+1)+(\max (a,b)+h)c(k+1)-a\cdot \min (b,c)\cdot (k+1)$.
\item Finally, the total number $\chi_3$ of matchings in
$[\begin{tikzpicture}[scale=0.25] \draw (0,0) [fill] circle (.1); \draw (1,0) [fill] circle (.1); \draw[thick] (1,0) arc (0:180:0.5); \end{tikzpicture}\, ,\mathbf{k}(\mathbf{h}(\mathbf{a};\mathbf{b});\mathbf{c})]$ having 3 small edges is immediate to compute,
and we get $\chi_3 =abc(h+1)(k+1)$.
\end{itemize}

Summing up the above contribution, we then find the desired closed expression for $|[\begin{tikzpicture}[scale=0.25] \draw (0,0) [fill] circle (.1); \draw (1,0) [fill] circle (.1); \draw[thick] (1,0) arc (0:180:0.5); \end{tikzpicture}\, ,\mathbf{k}(\mathbf{h}(\mathbf{a};\mathbf{b});\mathbf{c})]|$.

\bigskip

Our last example concerns a class of noncrossing matchings defined in a recursive fashion.
Before introducing them, we state an easy, but useful, lemma whose proof is left to the reader.

\begin{lemma}\label{equivalent}
Let $\sigma$ and $\tau$ be any matchings. Then the following are equivalent:
\begin{itemize}
\item[(i)] $\sigma \leq \tau$;
\item[(ii)] $\begin{tikzpicture}[scale=0.25, baseline=-0.3mm] \node[below] at (1.5,1) {$\sigma$}; \draw (0,0) [fill] circle (.1); \draw (3,0) [fill] circle (.1); \draw[thick] (3,0) arc (0:180:1.5); \end{tikzpicture} \leq \begin{tikzpicture}[scale=0.25, baseline=-0.3mm] \node[below] at (1.5,1) {$\tau$}; \draw (0,0) [fill] circle (.1); \draw (3,0) [fill] circle (.1); \draw[thick] (3,0) arc (0:180:1.5); \draw (4,0) [fill] circle (.1); \draw (6,0) [fill] circle (.1); \draw[thick] (6,0) arc (0:180:1); \end{tikzpicture}$
\item[(iii)] $\sigma \; \begin{tikzpicture}[scale=0.25] \draw (0,0) [fill] circle (.1); \draw (2,0) [fill] circle (.1); \draw[thick] (2,0) arc (0:180:1); \end{tikzpicture} \leq \begin{tikzpicture}[scale=0.25, baseline=-0.3mm] \node[below] at (1.5,1) {$\tau$}; \draw (0,0) [fill] circle (.1); \draw (3,0) [fill] circle (.1); \draw (5,0) [fill] circle (.1); \draw (6,0) [fill] circle (.1); \draw[thick] (5,0) arc (0:180:1); \draw[thick] (6,0) arc (0:180:3); \end{tikzpicture}$
\end{itemize}
\end{lemma}

Set $\tau_0 =\emptyset$. For every $n>0$, define:
\begin{itemize}
\item $\tau_{2n-1}=\tau_{2n-2} \; \begin{tikzpicture}[scale=0.25] \draw (0,0) [fill] circle (.1); \draw (2,0) [fill] circle (.1); \draw[thick] (2,0) arc (0:180:1); \end{tikzpicture}$ , and
\item $\tau_{2n}= \begin{tikzpicture}[scale=0.25, baseline=-0.3mm] \node[below] at (2.5,1) {$\tau_{2n-1}$}; \draw (0,0) [fill] circle (.1); \draw (5,0) [fill] circle (.1); \draw[thick] (5,0) arc (0:180:2.5); \end{tikzpicture}$ .
\end{itemize}

Denote with $f_{n}$ the cardinality of the interval
$[\begin{tikzpicture}[scale=0.25] \draw (0,0) [fill] circle (.1); \draw (1,0) [fill] circle (.1); \draw[thick] (1,0) arc (0:180:0.5); \end{tikzpicture}\, ,\tau_n ]$
and with $f_{n,k}$ the number of elements having $k$ edges of the same interval, for $k>0$. In particular, it is clear that $f_{n,k}=0$ for $n<k$ and whenever $n\leq 0$ or $k\leq 0$
(actually, when $n=k=0$, we set $f_{n,k}=0$ by convention).
In the next proposition we give closed formulas for such quantities.

\begin{proposition}
Let $n>0$ and $0<k\leq n$, and denote with $\varphi_n$ the $n$-th Fibonacci number. Then:
\begin{itemize}
\item[(i)] $f_{n,k}=\sum_{i=0}^{n-1}\binom{k-1}{n-k-1}$;
\item[(ii)] $f_n =\varphi_{n+2}-1$.
\end{itemize}
\end{proposition}

\proof  We use the following notations: $A_{n,k}$ is the set of all matchings in
$[\begin{tikzpicture}[scale=0.25] \draw (0,0) [fill] circle (.1); \draw (1,0) [fill] circle (.1); \draw[thick] (1,0) arc (0:180:0.5); \end{tikzpicture}\, ,\tau_{2n}]$ having $k$ edges,
$B_{n,k}$ is the set of all matchings in
$[\begin{tikzpicture}[scale=0.25] \draw (0,0) [fill] circle (.1); \draw (1,0) [fill] circle (.1); \draw[thick] (1,0) arc (0:180:0.5); \end{tikzpicture}\, ,\tau_{2n-1}]$ having $k$ edges,
and $C_{n,k}$ is the set of all matchings of the form
$\begin{tikzpicture}[scale=0.25, baseline=-0.3mm] \node[below] at (1.5,1) {$\sigma$}; \draw (0,0) [fill] circle (.1); \draw (3,0) [fill] circle (.1); \draw[thick] (3,0) arc (0:180:1.5); \end{tikzpicture}$ ,
with $\sigma \in B_{n,k-1}$.
We then have that $f_{2n,k}=|A_{n,k}|=|B_{n,k}|+|C_{n,k}|-|B_{n,k}\cap C_{n,k}|$.
By definition, we have $|B_{n,k}|=f_{2n-1,k}$ and clearly $|C_{n,k}|=f_{2n-1,k-1}$.
Furthermore, as a consequence of Lemma \ref{equivalent} and of the specific shape of the matchings under consideration,
the set $B_{n,k}\cap C_{n,k}$ is precisely the set of matchings of the form
$\begin{tikzpicture}[scale=0.25, baseline=-0.3mm] \node[below] at (1.5,1) {$\sigma$}; \draw (0,0) [fill] circle (.1); \draw (3,0) [fill] circle (.1); \draw[thick] (3,0) arc (0:180:1.5); \end{tikzpicture}$~,
with $\sigma \in B_{n-1,k-1}$, hence $|B_{n,k}\cap C_{n,k}|=f_{2n-3,k-1}$.
We thus get the recurrence relation $f_{2n,k}=f_{2n-1,k}+f_{2n-1,k-1}-f_{2n-3,k-1}$.
Using a completely similar argument, we can also prove the analogous recurrence $f_{2n-1,k}=f_{2n-2,k}+f_{2n-2,k-1}-f_{2n-4,k-1}$.
Summing up, we thus have the following recurrence relation, which holds for all $n,k\geq 2$:
\begin{equation}\label{gen_rec}
f_{n,k}=f_{n-1,k}+f_{n-1,k-1}-f_{n-3,k-1}.
\end{equation}

Together with the starting condition $f_{1,1}=1$,
formula (\ref{gen_rec}) allows us to compute the generating function $F(x,y)=\sum_{n,k\geq 0}f_{n,k}x^n y^k$.
Indeed, using standard arguments, our recurrence translates into the functional equation
$$
F(x,y)=xy+xF(x,y)+xyF(x,y)-x^3 F(x,y),
$$
which gives
$$
F(x,y)=\frac{xy}{1-x-xy+x^3y}.
$$

It turns out that $F(x,y)=xyG(x,y)$, where $G(x,y)$ is the generating function given in \cite{S} for the number triangle A004070:
from there, we deduce the desired closed form given in $(i)$ for $f_{n,k}$.
Moreover, denoting with $\Phi (x)=\sum_{n\geq 0}\varphi_n x^n$ the generating function of Fibonacci numbers,
it is easy to see that
$$
\Phi (x)-\frac{x}{1-x}=\frac{x}{1-x-x^2}-\frac{x}{1-x}=x^2 f(x,1),
$$
which proves $(ii)$.
\endproof

\section{Conclusion and further work} \label{CFW}

The enumerative combinatorics of the matching pattern poset remains still largely unknown.
Although some major efforts to enumerate pattern avoiding matchings have already been spent, as mentioned in Section $\ref{Work}$, the enumeration of most classes of matchings avoiding a single pattern of order three is still lacking.
To this regard, in the present paper we have introduced the notion of unlabeled pattern,
and we have enumerated matchings avoiding the unlabeled patterns [123123], [112323] and [123132], respectively.
However we did not succeed in finding a formula for the number of matchings avoiding the remaining two unlabeled patterns of order three,
namely [123321] and [112233], although matchings in the former class seem to have a rather neat combinatorial structure.

In Section $\ref{Int}$ we have started the investigation of the combinatorial structure of intervals in the matching pattern poset,
with special emphasis on enumerative issues. However, all important general questions concerning this topic are completely unanswered yet.
How many elements does a generic interval contain? How many (saturated) chain of fixed length? What is the M\"obius function?
In which cases an interval has a (possibly distributive) lattice structure?
Notice that the subposet of noncrossing matchings is isomorphic to the pattern order on 231-avoiding permutations
(this is rather easy to show, see also \cite{AB}).
This can be useful, for instance, in the computation of the M\"obius function, since the results developed in \cite{BJJS} can be applied.
However, it is possible (and maybe likely) that the specific combinatorial structure of matchings may help in finding neater formulas.

\end{document}